\documentclass{article}
\usepackage{amsmath}
\usepackage{amsfonts,amssymb}
\input amssym.def

\newcommand{\lie}{{\bf lie}}

\newcommand{\cocomm}{{\bf cocomm}}

\newcommand{\Coalg}{{\bf Coalg}\,}

\numberwithin{equation}{section}

\newcommand{\sgn}{{\rm s g n}}

\newcommand{\pr}{{\rm p r }}

\newcommand{\Hom}{{\rm H o m }}

\newcommand{\dia}{\diamond}

\newcommand{\tH}{{\widetilde{H}}}

\newcommand{\tal}{\widetilde{\alpha}}

\newcommand{\SM}{{\cal S}M}

\newcommand{\Omb}{\Om^{\bul}}

\newcommand{\OM}{\cO_M}

\newcommand{\Cbu}{C^{\bullet}}
\newcommand{\Cbd}{C_{\bullet}}

\newcommand{\Linf}{L_{\infty}}

\newcommand{\al}{{\alpha}}
\newcommand{\la}{{\lambda}}
\newcommand{\h}{{\hbar}}
\newcommand{\bul}{{\bullet}}

\newcommand{\mb}{{\mathfrak{b}}}

\newcommand{\Om}{{\Omega}}
\newcommand{\Si}{{\Sigma}}

\newcommand{\ga}{{\gamma}}

\newcommand{\ve}{{\varepsilon}}

\newcommand{\cF}{{\cal F}}
\newcommand{\pa}{{\partial}}

\newcommand{\cC}{{\cal C}}

\newcommand{\cL}{{\cal L}}

\newcommand{\cH}{{\cal H}}

\newcommand{\cU}{{\cal U}}

\newcommand{\cO}{{\cal O}}

\newcommand{\bbF}{{\Bbb F}}

\newcommand{\bbR}{{\Bbb R}}

\newcommand{\La}{{\Lambda}}

\newcommand{\D}{{\Delta}}

\newcommand{\tF}{{\widetilde{F}}}

\date{}
\newtheorem{defi}{Definition}

\newtheorem{lem}{Lemma}

\newtheorem{pred}{Proposition}

\title{Erratum to: ``A Proof of Tsygan's Formality Conjecture for an Arbitrary Smooth Manifold''}

\author{Vasiliy A. Dolgushev}

\begin{document}

\large

\maketitle

\begin{abstract} 
Boris Shoikhet noticed that the proof of 
lemma 1 in section 2.3 of \cite{thesis}
contains an error. In this note I give a correct proof of 
this lemma which was suggested to me by 
Dmitry Tamarkin. 
The correction does not change the results of \cite{thesis}.
\end{abstract}

\section{Introduction}
In this note I give a correct proof of lemma 1 
from section 2.3 in \cite{thesis}. This proof was 
kindly suggested to me by Dmitry Tamarkin and it 
is based on the interpretation of $\Linf$-morphisms 
as Maurer-Cartan elements of an auxiliary 
$\Linf$-algebra.  

The notion of partial homotopy proposed in 
section 2.3 in \cite{thesis} is poorly defined 
and this note should be used as a replacement of 
section 2.3 in \cite{thesis}. The main result of this
section (lemma 1) is used in section 5.2 of \cite{thesis}
in the proof of theorem 6. Since the statement of 
the lemma still holds so does the statement of 
theorem 6 as well as all other results of \cite{thesis}.
  
In section 2 of this note I recall the notion of 
an $\Linf$-algebra and the notion of a Maurer-Cartan 
element. In section 3, I give the interpretation of  
$\Linf$-morphisms as Maurer-Cartan elements of 
an auxiliary $\Linf$-algebra and use it to define 
homotopies between $\Linf$-morphisms. 
In section 4 I formulate and prove lemma 1 from 
section 2.3 of \cite{thesis}. Finally, in the 
concluding section, I give a model category 
interpretation of the homotopies between 
$\Linf$-morphisms. 

~\\
{\bf Notation.} I use the notation from 
\cite{thesis}. The underlying symmetric monoidal category
is the category of cochain complexes. For this
reason I sometimes omit the combination 
``DG'' (differential graded) talking about 
(co)operads and their (co)algebras. 
For a (co)operad $\cO$ I denote by $\bbF_{\cO}$
the corresponding Schur functor.  
$s\, K$ denotes the suspension of the complex $K$. 
In other words, 
$$
s \, K = s \otimes K\,, 
$$
where $s$ is the one-dimensional vector space 
placed in degree $+1$. Similarly, 
$$
s^{-1} K = s^{-1} \otimes K\,,
$$
where $s^{-1}$ is the one-dimensional vector space 
placed in degree $-1$\,.  $\lie$ is the operad of 
Lie algebras and 
$\cocomm$ is the cooperad of cocommutative coalgebras.

By ``suspension'' of a (co)operad $\cO$
I mean the (co)operad $\La(\cO)$
whose $m$-th space is
\begin{equation}
\label{susp-op}
\La(\cO)(m) = \Si^{1-m} \cO(m) \otimes \sgn_{m}\,,
\end{equation}
where $\sgn_{m}$ is the sign representation of
the symmetric group $S_m$\,.

~\\
{\bf Acknowledgment.}
I would like to thank E. Getzler, V. Hinich,  B. Shoikhet, 
and D. Tamarkin for useful discussions. In particular, 
it is D. Tamarkin who suggested to me the
correct proof of lemma \ref{styag} and it is V. Hinich
who explained to me the model category interpretation 
of homotopies between $\Linf$-morphisms.

\section{$\Linf$-algebras and Maurer-Cartan elements}
Let me recall from \cite{GK} that an $\Linf$-algebra 
structure on a graded vector space $\cL$ is a 
degree $1$ codifferential $Q$ on the 
colagebra
$\bbF_{\La\cocomm}(\cL)$ cogenerated by $\cL$\,.
Following \cite{thesis} I denote 
the DG coalgebra $(\bbF_{\La\cocomm}(\cL),Q)$ by 
$C(\cL)$: 
\begin{equation}
\label{C-cL}
C(\cL) = (\bbF_{\La\cocomm}(\cL),Q)\,.
\end{equation}

A morphism $F$ from an $\Linf$-algebra 
$(\cL,Q)$ to an $\Linf$-algebra $(\cL^{\dia}, Q^{\dia})$
is by definition a morphism of (DG) coalgebras
\begin{equation}
\label{FFF}
F :  C(\cL) \to C(\cL^{\dia})\,.
\end{equation}

Since
$$
\bbF_{\La\cocomm}(\cL) = s\, \bbF_{\cocomm}(s^{-1} \cL)
$$
the vector space of $C(\cL)$ can 
be identified with the exterior algebra 
$\wedge^{\bul}\cL$ and for a graded vector 
space $V$ a map 
$$
f : \bbF_{\La\cocomm}(\cL) \to V 
$$
of degree $|f|$
can be identified with the infinite collection
of maps 
$$
f_n : \cL^{\otimes\, n} \to V\,, \qquad n\ge 1\,,
$$
where each map $f_n$ has degree $|f|+ 1 - n$ and 
$$
f_n (\dots, \ga, \ga', \dots) = 
- (-1)^{|\ga| |\ga'|} f_n(\dots, \ga', \ga, \dots)
$$ 
for every pair of elements $\ga, \ga' \in \cL$\,.

Due to proposition 2.14 in \cite{GJ} every coderivation 
of $\bbF_{\La\cocomm}(\cL)$ is uniquely determined 
by its composition with the projection 
\begin{equation}
\label{pr}
\pr_{\cL} : \bbF_{\La\cocomm}(\cL) \to \cL
\end{equation}
from $\bbF_{\La\cocomm}(\cL)$ onto cogenerators.

In particular, the codifferential $Q$ of the 
coalgebra $C(\cL)$ is uniquely determined by 
the infinite collection of maps 
\begin{equation}
\label{Q-n}
Q_n = \pr_{\cL}\circ Q \Big|_{\wedge^n \cL} 
:  \wedge^n \cL \to \cL\,, 
\end{equation}
such that $Q_n$ has degree $2 - n$\,.
In \cite{thesis} $Q_n$ are called structure maps 
of the $\Linf$-algebra $\cL$\,.

The equation $Q^2=0$ is equivalent to an infinite 
collection of quadratic equations on the maps 
(\ref{Q-n}). The precise form of these equations 
can be found in definition $4.1$ in \cite{Ezra}. 

One of the obvious equations implies that the 
structure map of the first level
$Q_1$ is a degree $1$ differential of $\cL$\,. 
Thus an $\Linf$-algebra can be thought of 
as an algebra over an operad in the category 
of cochain complexes. 

Every $\Linf$-algebra $\cL$ is equipped with a 
natural decreasing filtration:
$$
\cL = \cF^1_{lc} \cL \supset  \cF^2_{lc} \cL 
\supset \cF^3_{lc} \cL \supset \dots 
$$
$$
\cF^i_{lc} = \bigoplus_{i_1+ \dots + i_k = i}
Q_k( \cF^{i_1}_{lc} \cL, \cF^{i_2}_{lc} \cL, \dots, 
\cF^{i_k}_{lc} \cL )\,,
$$
which is called the lower central filtration.  

\begin{defi}[E. Getzler \cite{Ezra}]
\label{nilp}
An $\Linf$-algebra $\cL$ is nilpotent if the 
lower central filtration of $\cL$ terminates, that is, 
if $\cF^i_{lc} \cL =0$ for $i \gg 0$\,.
\end{defi}
Furthermore, an $\Linf$-algebra is called pronilpotent 
if it is a projective limit of nilpotent 
$\Linf$-algebras.

For a (pro)nilpotent  $\Linf$-algebra $\cL$
it makes sense to speak about its Maurer-Cartan 
elements: 
\begin{defi}[E. Getzler \cite{Ezra}]
A Maurer-Cartan $\pi$ of a pronilpotent 
$\Linf$-algebra $(\cL, Q)$ is a degree $1$ element 
of $\cL$ satisfying the equation  
\begin{equation}
\label{MC}
\sum_{n=1}^{\infty} \frac{1}{n!} Q_n (\pi, \pi, \dots, \pi) = 0\,.
\end{equation}
\end{defi}
Let me remark that the infinite sum 
in (\ref{MC}) is well defined since $\cL$ is 
pronilpotent. 

Every Maurer-Cartan element $\pi$ of $\cL$ can be 
used to modify the $\Linf$-algebra structure 
on $\cL$\,. This modified structure is called 
the $\Linf$-structure twisted by the Maurer-Cartan 
$\pi$ and its structure maps are given by
\begin{equation}
\label{Q-pi}
Q^{\pi}_n (\ga_1, \dots, \ga_n) =
\sum_{m=1}^{\infty} \frac{1}{m!}
Q_{m+n} (\pi, \dots , \pi, \ga_1, \dots, \ga_n)\,, 
\qquad \ga_i \in \cL\,. 
\end{equation}
It is equation (\ref{MC}) which implies that 
the maps (\ref{Q-pi}) define an $\Linf$-algebra 
structure on $\cL$\,.

Two Maurer-Cartan elements $\pi_0$  
and $\pi_1$ are called equivalent if there is an element 
$\xi\in \cL^0$ such that the solution of the 
equation 
\begin{equation}
\label{MC-homot}
\frac{d}{d t} \pi_t = Q_1^{\pi_t} (\xi)
\end{equation}
connects $\pi_0$ and $\pi_1$:
$$
\pi_t \Big|_{t=0} = \pi_0\,, 
\qquad 
\pi_t \Big|_{t=1} = \pi_1\,.
$$

\section{$\Linf$-morphisms and their homotopies}

I will need the following 
auxiliary statement: 
\begin{pred}
\label{Hom-s-cocomm}
Let $\cO$ be an operad and $A$ be an algebra
over $\cO$. If $B$ is a (DG) cocommutative coalgebra
then the cochain complex
\begin{equation}
\label{H-B-A}
\cH_{B,A}= \Hom(B,A)
\end{equation}
of all linear maps from $B$ to $A$
has a natural structure of an algebra over $\cO$\,. 
\end{pred}
{\bf Proof.} The $\cO$-algebra structure on 
$A$ is by definition the map (of complexes)
\begin{equation}
\label{mu-cO-A}
\mu_A : \bbF_{\cO}(A) \to A  
\end{equation}
making the following diagrams commutative: 
\begin{equation}
\label{mu-mu-cO}
\begin{array}{ccc}
\bbF_{\cO} ( \bbF_{\cO}(A)) &\,
\stackrel{\bbF_{\cO}(\mu_A)}{\longrightarrow}\, &
 \bbF_{\cO}(A) \\[0.3cm]
\downarrow^{~\mu_{\cO}(A)} & ~ & \downarrow^{~\mu_A} \\[0.3cm]
  \bbF_{\cO}(A) &\, \stackrel{\mu_A}{\longrightarrow} \, & A\,,  
\end{array}
\end{equation}
\begin{equation}
\label{mu-u-cO}
\begin{array}{ccc}
A  & \, \stackrel{u_{\cO}(A)}{\longrightarrow}\, &  
\bbF_{\cO}(A) \\[0.3cm]
~ & \searrow^{~ {\rm i d}}  & ~~\downarrow^{~\mu_A} \\[0.3cm]
 ~  &  ~  &  A
\end{array}
\end{equation}
where $\mu_{\cO}$ and $u_{\cO}$ are
the transformation of functors
$$
\mu_{\cO} : \bbF_{\cO}\circ \bbF_{\cO} \to \bbF_{\cO}\,,
$$
$$
u_{\cO} : {\rm I d} \to \bbF_{\cO}
$$ 
defined by the operad structure on $\cO$\,.
The map $\mu_A$ is called the multiplication.

For every $n>1$
the comultiplication $\D$ in $B$ provides me 
with the following map 
$$
\D^{(n)} : B \to B^{\otimes \, n}
$$
\begin{equation}
\label{Delta-n}
\D^{(n)} X = 
(\D \otimes 1^{\otimes \, (n-2)}) \dots
(\D \otimes 1 \otimes 1) (\D \otimes 1)\D \, X
\end{equation}
Using this map and the $\cO$-algebra structure 
on $A$\,, I define the $\cO$-algebra structure on 
$\cH_{B,A}$ (\ref{H-B-A}) by 
\begin{equation}
\label{mu-Hom}
\mu(v, \ga_1, \dots, \ga_n; X) = 
\mu_A (v) [\ga_1\otimes \dots \otimes \ga_n \,(\D^{(n)} X ) ]\,, 
\end{equation}
where $v\in \cO(n)$\,, $\ga_i\in \Hom(B,A)$\,,
and $X\in B$\,.

The equivariance with respect to the action of the 
symmetric group follows from the cocommutativity 
of the comultiplication on $B$\,.

The commutativity of the diagram 
\begin{equation}
\label{mu-mu-cO-Hom}
\begin{array}{ccc}
\bbF_{\cO} ( \bbF_{\cO}(\cH_{B,A})) &\,
\stackrel{\bbF_{\cO}(\mu)}{\longrightarrow}\, &
 \bbF_{\cO}(\cH_{B,A}) \\[0.3cm]
\downarrow_{~\mu_{\cO}(\cH_{B,A})} & ~ & ~~\downarrow^{~\mu} \\[0.3cm]
\bbF_{\cO}(\cH_{B,A}) &\, \stackrel{\mu}{\longrightarrow} \, & \cH_{B,A}\,,  
\end{array}
\end{equation}
follows from the commutativity of (\ref{mu-mu-cO})
and the associativity of the comultiplication 
in $B$\,. 

The commutativity of the diagram
\begin{equation}
\label{mu-u-cO-Hom}
\begin{array}{ccc}
\cH_{B,A}  & \, \stackrel{u_{\cO}(\cH_{B,A})}{\longrightarrow}\, &  
\bbF_{\cO}(\cH_{B,A}) \\[0.3cm]
~ & \searrow^{~ {\rm i d}}  & ~~\downarrow^{~\mu} \\[0.3cm]
 ~  &  ~  &  \cH_{B,A}
\end{array}
\end{equation}
and the compatibility of $\mu$ (\ref{mu-Hom}) 
with the differential are obvious. $\Box$

Since
$$
\bbF_{\La\cocomm}(\cL) =
s\, \bbF_{\cocomm} (s^{-1} \cL) 
$$
for every $\Linf$-algebra $\cL^{\dia}$
proposition \ref{Hom-s-cocomm}
gives me a $\Linf$-structure on the 
cochain complex
\begin{equation}
\label{ona}
\cU = s\, \Hom(C(\cL), \cL^{\dia})\,.
\end{equation}

This algebra $\cU$ can be equipped with the 
following decreasing filtration: 
$$
\cU = \cF^1 \cU \supset \cF^2 \cU \supset \dots
\supset \cF^k \cU  \supset \dots 
$$
\begin{equation}
\label{filtr}
\cF^k \cU = \{f \in \Hom(\wedge^\bul \cL, \cL^{\dia})
\quad | \quad f\Big |_{\wedge^{< k} \cL} = 0 \}\,. 
\end{equation}
It is not hard to see that this filtration is 
compatible with the $\Linf$-algebra structure on $\cU$\,.
Furthermore, since $\cU = \cF^1 \cU$, for every $k$
the $\Linf$-algebra $\cU / \cF^k \cU$ is nilpotent.   
On the other hand, 
\begin{equation}
\label{polna}
\cU = \lim_{k} \cU / \cF^k \cU\,,
\end{equation}
and hence, the $\Linf$-algebra $\cU$
is pronilpotent and the notion of a Maurer-Cartan 
element of $\cU$ makes sense.

My next purpose is to identify the Maurer-Cartan elements 
of the $\Linf$-algebra $\cU$ (\ref{ona}) with $\Linf$-morphisms 
from $\cL$ to $\cL^{\dia}$ :
\begin{pred}
\label{Linf-morph}
$\Linf$-morphisms from $\cL$ to $\cL^{\dia}$
are identified with Maurer-Cartan elements 
of the $\Linf$-algebra $\cU$ (\ref{ona}) 
\end{pred}
{\bf Proof.} 
Since $ C(\cL^{\dia})$ is a cofree coalgebra,  
the map $F$ (\ref{FFF}) is uniquely determined by 
its composition $\pr_{\cL^{\dia}} \circ F$ with 
the projection $\pr_{\cL^{\dia}}$ (\ref{pr}).
This composition is a degree zero element of 
$\Hom(C(\cL), \cL^{\dia})$\,. Thus, since $\cU$ (\ref{ona}) 
is obtained from $\Hom(C(\cL), \cL^{\dia})$
by the suspension, every morphism $F$ (\ref{FFF}) is 
identified with a degree $1$ element of $\cU$\,. 
 
It remains to prove that the compatibility condition
\begin{equation}
\label{QF-FQ}
Q^{\dia} F = F Q
\end{equation}
of $F$ with the codifferentials $Q$ and $Q^{\dia}$
on $C(\cL)$ and $C(\cL^{\dia})$, respectively, is 
equivalent to the Maurer-Cartan equation (\ref{MC})
on $\pr_{\cL^{\dia}} \circ F$ viewed as an element 
of $\cU$\,.

It is not hard to see that 
\begin{equation}
\label{QF-FQ1}
\pr_{\cL^{\dia}} \circ (Q^{\dia} F - F Q) = 0\,.
\end{equation}
is equivalent to the Maurer-Cartan equation 
on the composition $\pr_{\cL^{\dia}} \circ F $
viewed as an element of $\cU$ (\ref{ona}). 

Thus, I have to show that
equation (\ref{QF-FQ1}) is equivalent to 
the compatibility condition (\ref{QF-FQ}).

For this, I denote by $\Psi$ the difference: 
$$
\Psi = Q^{\dia} F - F Q 
$$
and remark that 
\begin{equation}
\label{Psi-Delta}
\D \Psi = - (\Psi \otimes F + F \otimes \Psi ) \D\,,
\end{equation}
where $\D$ denotes the coproduct both in $C(\cL)$ 
and $C(\cL^{\dia})$\,.

The latter follows from the fact that $Q$ and $Q^{\dia}$
are coderivations and $F$ is a morphism of 
cocommutative coalgebras. 

Given a cooperad $\cC$, a pair of cochain
complexes $V$, $W$, a degree zero map   
$$
f: V \to W 
$$ 
and an arbitrary map 
$$
b : V  \to W
$$
I denote by $\pa(b,f)$ the following 
map\footnote{A similar construction was introduced at the beginning 
of section 2.2 in \cite{GJ}.} 
$$
\pa(b,f) : \bbF_{\cC}(V) \to \bbF_{\cC}(W) 
$$
\begin{equation}
\label{pa}
\begin{array}{c}
\pa(b,f) (\ga, v_1, v_2, \dots, v_n)  = \\[0.3cm] 
\displaystyle
\sum_{i=1}^n 
(-1)^{|b|(|\ga|+ |v_1| + \dots + |v_{i-1}|)}
(\ga, f(v_1), \dots, f(v_{i-1}), 
b(v_i), f(v_{i+1}), \dots, f(v_n))\,,
\end{array}
\end{equation}
$$
\ga\in \cC(n)\,, \qquad v_i \in V\,, 
$$
where $|\ga|$, $|b|$, $|v_j|$ are, respectively, degrees 
of $\ga$, $b$, and $v_j$\,. The equivariance of (\ref{pa}) 
with respect to permutations is obvious. 

It is not hard to see that condition 
(\ref{Psi-Delta}) is equivalent to commutativity 
of the following diagram 
\begin{equation}
\label{nu-Psi-F}
\begin{array}{ccc}
\bbF_{\La\cocomm}(\cL) &\,
\stackrel{\Psi}{\longrightarrow}\, &
 \bbF_{\La\cocomm}(\cL^{\dia}) \\[0.3cm]
\downarrow^{~\nu} & ~ & \downarrow^{~\nu} \\[0.3cm]
\bbF_{\La\cocomm}
(\bbF_{\La\cocomm}(\cL)) &\, 
\stackrel{\pa(\Psi, F)}{\longrightarrow} \, & 
\bbF_{\La\cocomm}
(\bbF_{\La\cocomm}(\cL^{\dia}))\,,  
\end{array}
\end{equation}
where $\nu$ is the coproduct of the 
cotriple $\bbF_{\La\cocomm}$\,.

Since the functor $\bbF_{\La\cocomm}$ with the transformations
$\nu:\bbF_{\La\cocomm} \to \bbF_{\La\cocomm} \circ \bbF_{\La\cocomm} $
and $\pr : \bbF_{\La\cocomm} \to {\rm Id} $ form a 
cotriple\footnote{See, for example, section 1.7 in \cite{GJ}.}, 
the following diagram   
\begin{equation}
\label{pr-Psi-F}
\begin{array}{ccc}
\bbF_{\La\cocomm}(\cL^{\dia}) & ~ & ~  \\[0.3cm]
\downarrow^{~\nu} & \searrow^{~{\rm i d}} & ~ \\[0.3cm]
\bbF_{\La\cocomm}
(\bbF_{\La\cocomm}(\cL^{\dia})) 
&\, \stackrel{p}{\longrightarrow}\, 
& \bbF_{\La\cocomm}(\cL^{\dia})\,,
\end{array}
\end{equation}
with $p$ being $\bbF_{\La\cocomm}(\pr_{\cL^{\dia}})$, commutes.

Attaching this diagram to (\ref{nu-Psi-F}) 
I get the commutative diagram
\begin{equation}
\label{xorom}
\begin{array}{ccccc}
\bbF_{\La\cocomm}(\cL) &
\stackrel{\Psi}{\longrightarrow} &
 \bbF_{\La\cocomm}(\cL^{\dia}) & ~ & ~ \\[0.3cm]
\downarrow^{~\nu} & ~ & \downarrow^{~\nu} & \searrow^{~{\rm i d}} 
& ~   \\[0.3cm]
\bbF_{\La\cocomm}
(\bbF_{\La\cocomm}(\cL)) & 
\stackrel{\pa(\Psi, F)}{\longrightarrow}  & 
\bbF_{\La\cocomm}
(\bbF_{\La\cocomm}(\cL^{\dia})) & 
  \stackrel{p}{\longrightarrow} & 
 \bbF_{\La\cocomm}(\cL^{\dia})\,,
\end{array}
\end{equation}
where, as above, $p = \bbF_{\La\cocomm}(\pr_{\cL^{\dia}})$\,.

Hence, 
$$
\Psi = \bbF_{\La\cocomm}(\pr_{\cL^{\dia}}) 
\circ \pa(\Psi, F) \circ \nu \,.
$$
On the other hand 
$$
\bbF_{\La\cocomm}(\pr_{\cL^{\dia}}) \circ \pa(\Psi, F) = 
\pa (\pr_{\cL^{\dia}} \circ \Psi, \pr_{\cL^{\dia}} \circ F)\,.
$$
Therefore,
$$
\Psi = \pa (\pr_{\cL^{\dia}} \circ \Psi, 
\pr_{\cL^{\dia}} \circ F) \circ \nu
$$
and $\Psi$ vanishes if and only if so 
does the composition $\pr_{\cL^{\dia}} \circ \Psi$\,.

This concludes the proof of the proposition. $\Box$

The identification proposed in the above proposition 
allows me to introduce a notion of homotopy between 
two $\Linf$-morphisms. Namely,
\begin{defi}
\label{homot}
$\Linf$-morphisms $F$ and $\tF$ from $\cL$ to $\cL^{\dia}$
are called homotopic if the corresponding Maurer-Cartan 
elements of the $\Linf$-algebra $\cU$ (\ref{ona}) are 
equivalent. 
\end{defi}

\section{Lemma $1$ from \cite{thesis}}

Let me denote by $F_n$ the components 
$$
F_n  : \wedge^n \, \cL \to \cL^{\dia} 
$$ 
\begin{equation}
\label{F-n}
F_n = \pr_{\cL^{\dia}} \circ F \Big |_{\wedge^n\, \cL}
\end{equation}
of the composition $\pr_{\cL^{\dia}} \circ F$\,, 
where $\pr_{\cL^{\dia}}$ is the projection from 
$\bbF_{\La\cocomm}(\cL^{\dia})$ onto cogenerators. 
In \cite{thesis} the maps (\ref{F-n}) are 
called structure maps of the $\Linf$-morphism 
(\ref{FFF}). 

The compatibility condition (\ref{QF-FQ}) implies that 
the structure map $F_1$ of the first level is morphism 
of complexes: 
$$
F_1 : \cL \to \cL^{\dia} \,, 
\qquad 
Q^{\dia}_1 F_1 = F_1 Q_1\,. 
$$
By definition, an $\Linf$-morphism $F$ is a 
$\Linf$-quasi-isomorphism if the map $F_1$ is a 
quasi-isomorphism of the corresponding complexes.  

I can now prove the following lemma:
\begin{lem}
\label{styag}
Let
$$
F : C(\cL) \mapsto C(\cL^{\dia})
$$
be a quasi-isomorphism from
an $\Linf$-algebra $(\cL, Q)$ to
an $\Linf$-algebra $(\cL^{\dia}, Q^{\dia})$.
For $n\ge 1$ and any map
\begin{equation}
\label{tH}
\tH : \wedge^{n}\cL \mapsto \cL^{\dia}
\end{equation}
of degree $-n$ one can construct
a quasi-isomorphism
$$
\tF : C(\cL) \mapsto C(\cL^{\dia})
$$
such that for any $m<n$
\begin{equation}
\label{kak}
\tF_m = F_m \, :\, \wedge^m \cL\mapsto \cL^{\dia}
\end{equation}
and
$$
\tF_n (\ga_1, \dots, \ga_n)= F_n (\ga_1, \dots, \ga_n) +
$$
\begin{equation}
\label{kak1}
Q^{\dia}_1 \tH(\ga_1, \dots, \ga_n) -
(-)^{n} \tH (Q_1(\ga_1), \ga_2, \dots, \ga_n) - \dots
\end{equation}
$$
\dots -(-)^{n+k_1+ \dots + k_{n-1}}
\tH(\ga_1, \dots, \ga_{n-1}, Q_1(\ga_n))\,,
$$
where $\ga_i \in \cL^{k_i}$\,.
\end{lem}
{\bf Proof.} Let $Q^{\cU}$ denote the $\Linf$-algebra 
structure on $\cU$ (\ref{ona}). Let $\al$  
be the Maurer-Cartan elements of $\cU$
corresponding to the $\Linf$-morphism $F$. 

By setting 
\begin{equation}
\label{xi}
\xi \Big|_{\wedge^m \cL} = 
\begin{cases}
\tH \,, ~ {\rm if} ~ m = n\,, \\
0\,, ~ {\rm otherwise}
\end{cases}
\end{equation}
I define an element $\xi\in \cU$ of degree $0$. 
By definition of the filtration (\ref{filtr}) 
the element $\xi$ belongs to $\cF^n \cU$

Let $\al_t$ be the unique path of Maurer-Cartan 
elements defined by 
\begin{equation}
\label{al-t}
\frac{d}{d t} \al_t = (Q^{\cU})^{\al_t}_1 (\xi)\,, 
\qquad \al_t \Big|_{t = 0} = \al\,. 
\end{equation}

The unique solution $\al_t$ of (\ref{al-t}) can 
be found by iterating the following equation in 
degrees in $t$
\begin{equation}
\label{al-t1}
\al_t =  \al + \int_0^t (Q^{\cU})^{\al_{\tau}}_1 (\xi) d \tau\,. 
\end{equation}
Since the $\Linf$-algebra $\cU$ is pronilpotent 
the recurrent procedure (\ref{al-t1}) converges.

It is not hard to see that, since $\xi \in \cF^n \cU$, 
\begin{equation}
\label{al-t-n}
\al_t -  \al \in \cF^{n}\cU    
\end{equation}
and
\begin{equation}
\label{al-t-n1}
\al_t -  (\al + t Q^{\cU}_1 (\xi)) \in \cF^{n+1} \cU\,.   
\end{equation}

Let $\tF$ be the $\Linf$-morphism from $\cL$ to 
$\cL^{\dia}$ corresponding to the Maurer-Cartan 
element 
$$
\tal = \al_t \Big|_{t=1}\,.
$$

Equation (\ref{al-t-n}) implies (\ref{kak}) and
equation (\ref{al-t-n1}) implies (\ref{kak1})\,.
It is obvious that,
since $F$ is a quasi-isomorphism, so is $\tF$\,.  

The lemma is proved.  $\Box$
 
\section{Model category interpretation of the 
homotopies} 
In \cite{Hinich} V. Hinich showed that the category 
$\Coalg$ of 
unital (unbounded) DG cocommutative coalgebras can 
be equipped with a structure of the closed model 
category. The definition is based on Quillen's  
functor $\Om_{\La\lie}$ \cite{Q} from the category $\Coalg$ 
to the category 
of DG algebras over the operad $\La\lie$\,. Namely, the 
cofibrations in $\Coalg$ are injective maps and 
weak equivalences are maps $f$ such that $\Om_{\La\lie}(f)$
is a quasi-isomorphism.

In this section I give an interpretation\footnote{This interpretation 
was explained to me by V. Hinich.} 
of the homotopies 
between morphisms of $\Linf$-algebras $(\cL, Q)$ and 
$(\cL^{\dia}, Q^{\dia})$ in terms of this model 
category structure.  

First, I notice that for every $\Linf$-algebra 
the coalgebra $C(\cL)$ can be easily upgraded 
to a unital coalgebra in the sense of 
definition 2.1.1. in \cite{Hinich} by attaching 
the group-like element $u$ 
\begin{equation}
\label{attach}
C^+(\cL) = C(\cL) \oplus \bbR u
\end{equation}
with the properties
$$
Q (u) = 0\,, \qquad
\D u = u \otimes u\,, \qquad 
\ve(u) = 1\,, 
$$
where $\D$ is the comultiplication and 
$\ve$ is the counit.  

Second, one can similarly extend an
$\Linf$-morphism $F$ from $\cL$ to 
$\cL^{\dia}$ to a morphism $F^+$ between 
the corresponding unital coalgebras $C^+(\cL)$ 
and $C^+(\cL^{\dia})$\,. It is also obvious that 
every morphism of the unital coalgebras 
$C^+(\cL)$ and $C^+(\cL^{\dia})$ can be obtained in this way.
Thus I get a fully faithful embedding of the category of 
$\Linf$-algebras to the category of unital 
(unbounded) DG cocommutative coalgebras. 
I denote this embedding by $C^+$\,. 

I claim that
\begin{pred}
\label{CMC}
Two $\Linf$-morphisms $F$ and $\tF$ are
homotopic in the sense of definition \ref{homot}
if and only if $F^+$ and $\tF^+$ are homotopic 
in the closed model category  of unital 
(unbounded) DG cocommutative coalgebras.
\end{pred}
{\bf Proof.} It is obvious that for every 
$\Linf$-algebra the coalgebra $C^+(\cL)$ is a
cofibrant object in $\Coalg$\,. Furthermore, 
using the fact that the coalgebra $C^+(\cL)$ is 
free it is not hard to show that $C^+(\cL)$ is 
also fibrant.  Therefore it suffices to show that the 
homotopy between $F$ and $\tF$ given in definition
\ref{homot} is equivalent to the right homotopy 
between $F^+$ and $\tF^+$
with a fixed very good path object $C^+(\cL^{\dia})^I$ 
for $C^+(\cL^{\dia})$\,.  

I choose the following path object for $C^+(\cL^{\dia})$:
\begin{equation}
\label{path}
C^+(\cL^{\dia})^I = C^+(\cL^{\dia}\otimes \Omb(\bbR))\,, 
\end{equation}
where 
$$
\Omb(\bbR) = \bbR[t] \oplus \bbR[t] d t
$$ 
is the polynomial De Rham algebra of the real line
with $d t$ having degree $1$\,.

The natural embedding
\begin{equation}
\label{iota}
\iota : \cL^{\dia} \hookrightarrow \cL^{\dia}\otimes \Omb(\bbR)
\end{equation}
and the natural projections
$$
p_0 (X) = X \Big|_{t=0,\, d t =0} :   
\cL^{\dia}\otimes \Omb(\bbR)  \to \cL^{\dia}\,, 
$$
$$
p_1 (X) = X \Big|_{t=1,\, d t =0} :   
\cL^{\dia}\otimes \Omb(\bbR)  \to \cL^{\dia}\,, 
$$ 
provide me with the required morphisms
\begin{equation}
\label{C-iota}
C^+(\iota) : C^+(\cL^{\dia})  \hookrightarrow
C^+(\cL^{\dia}\otimes \Omb(\bbR))
\end{equation}
and 
$$
C^+(p_0) :
C^+(\cL^{\dia}\otimes \Omb(\bbR) )  \to 
C^+(\cL^{\dia})\,, 
$$
$$
C^+(p_1):
C^+(\cL^{\dia}\otimes \Omb(\bbR) )  \to 
C^+(\cL^{\dia})\,. 
$$ 

Since $\iota$ (\ref{iota}) is a quasi-isomorphism 
of $\Linf$-algebras $C^+(\iota)$ is weak equivalence.
Furthermore, $C^+(\iota)$ is also an embedding and 
hence is an acyclic cofibration. Thus (\ref{path}) 
is a very good path object.

Let $H^+$ be a morphism from $C^+(\cL)$ to 
$C^+(\cL^{\dia}\otimes \Omb(\bbR))$. Then $H^+$
is defined by the corresponding morphism 
\begin{equation}
\label{H}
H : C(\cL) \to C(\cL^{\dia}\otimes \Omb(\bbR))
\end{equation}
which is, in turn, uniquely determined by 
its composition 
$pr_{\cL^{\dia}\otimes \Omb(\bbR)}\circ H$
with the projection $pr_{\cL^{\dia}\otimes \Omb(\bbR)}$ 
from $C(\cL^{\dia}\otimes \Omb(\bbR))$ onto 
$\cL^{\dia}\otimes \Omb(\bbR)$\,. I denote this composition 
by $h$
\begin{equation}
\label{h}
h = pr_{\cL^{\dia}\otimes \Omb(\bbR)} \circ H : 
C(\cL) \to \cL^{\dia}\otimes \Omb(\bbR)\,.
\end{equation}

According to proposition \ref{Linf-morph} the 
element $h$ is a Maurer-Cartan element of the 
$\Linf$-algebra $s \, \Hom(C(\cL), \cL^{\dia}\otimes \Omb(\bbR))$\,.

Let me decompose the element $h$ as 
$$
h = h^0 + h^1 d t\,,
$$   
where $h^0$ and $h^1$ are elements in 
$\Hom(C(\cL), \cL^{\dia}[t])$ of 
degrees $0$ and $-1$, respectively. 
In terms of this decomposition, the Maurer-Cartan 
equation (\ref{MC}) for $h$ (\ref{h}) boils down to 
\begin{equation}
\label{MC1}
\sum_{n=1}^{\infty} \frac{1}{n!} 
Q^{\Hom}_n (h^0, \dots, h^0) = 0\,,
\end{equation}
\begin{equation}
\label{MC11}
\frac{\pa h^0}{\pa t}  = 
\sum_{n=1}^{\infty} \frac{1}{n!} 
Q^{\Hom}_{n+1} (h^0, \dots, h^0, h^1)\,,
\end{equation}
where $Q^{\Hom}_k$ are structure maps of the 
$\Linf$-algebra $s\, \Hom(C(\cL), \cL^{\dia})$

Equation (\ref{MC1}) tells me that $h^0$ is  
a Maurer-Cartan element of the $\Linf$-algebra 
$$
s\, \Hom(C(\cL), \cL^{\dia}[t])
$$ 
and $h^1$ defines an equivalence between the 
Maurer-Cartan element 
$$
h^0 \Big|_{t=0}  \in s\, \Hom(C(\cL), \cL^{\dia})\,, 
$$
and 
$$
h^0 \Big|_{t=1} \in s\, \Hom(C(\cL), \cL^{\dia})\,. 
$$
 
This consideration allows me to conclude 
that two $\Linf$-morphisms
$F$ and $\tF$ between the $\Linf$-algebras $\cL$ and
$\cL^{\dia}$ are homotopic in the sense of 
definition \ref{homot} if and only if 
there is a morphism of DG coalgebras
$$
H^+ : C^+(\cL) \to C^+(\cL^{\dia}\otimes \Omb(\bbR)) 
$$
such that 
$$
C^{+}(p_0)\circ H^+ = F^+ \,,
$$
and 
$$
C^{+}(p_1)\circ H^+ = \tF^+\,.
$$

Thus the proposition follows. $\Box$

\noindent\textsc{Department of Mathematics,
Northwestern University, \\
Evanston, IL 60208, USA \\
\emph{E-mail address:} {\bf vald@math.northwestern.edu}}

\end{document}